\documentclass[a4paper,11pt]{amsart}
\usepackage{amssymb,amsfonts,amsmath}
\usepackage{layout}
\usepackage[latin1]{inputenc}
\usepackage{amsthm}
\usepackage[bookmarksnumbered, colorlinks, plainpages]{hyperref}

\def\c2{\mathbb{C}^2}

\def\1{\bold{1}}

\theoremstyle{definition}

\theoremstyle{plain}
\newtheorem{def/not}{Definition/Notations}
\numberwithin{equation}{section}

\begin{document}

\title[Representations of the generalized operator entropy]{Continued fraction representations of the generalized operator entropy}
\author[S. Ahallal, S. Mennou, A. Kacha]{ Sarra Ahallal, Said Mennou  and Ali Kacha$^*$ }
\address{ Ibn Toufail University, Science faculty, Laboratory EDPAGS, Kenitra  $14000,$ Morocco.}
\address{$^*$ Corresponding author.}

{\bigskip}

\email{sarra.ahallal@uit.ac.ma}
\email{saidmennou@yahoo.fr}
\email{ali.kacha@uit.ac.ma}

\maketitle

\noindent {\bf { Abstract.\ }}
The direct calculation of the Generalized operator entropy  proves difficult by the
appearance of rational exponents of matrices. The main motivation of this
work is to overcome these difficulties and to present a practical and
efficient method for this calculation using its representation by the matrix continued fraction. At the end of our paper,
we deduce a continued fraction expansion of the Bregman operator divergence. Some numerical examples illutrating the theoretical result are discussed.

\rule[-10pt]{0pt}{10pt}\newline


{\bigskip}

\noindent
\large{\textbf{2. Introduction and motivation}}

{\bigskip}

\noindent
It is known that the real continued fraction expansions have the advantage that they converge
more rapidly than other numerical algorithms $[7,8].$ So the extension of continued fractions theory from real numbers to
the matrix case has seen several development and interesting applications, $[5].$
\\
\noindent The theory of operator means for positive and bounded linear operators on
a Hilbert space was initiated by T. Ando $[1]$
and established by him and F. Kubo in connection with Loweners theory for
the operator monotone functions. It is started from the presence
of the notion of parallel sum as a tool for analyzing multi-port electrical
networks in engineering, see $[2,3].$
\\
\noindent In $1850,$ Clausius, introduced the notion of entropy in
thermodynamics. Since then several extensions and reformulations
have been developed in various disciplines $[10,13].$ There
have been investigated the so called entropy inequalities by some
mathematicians, see $[14,15]$ and references therein.
\\
\\
\noindent For a real number $q \in \mathbb{R}, 0<q<1$ we define the path
$$
f_q(A,B)=A^{1/2}
(A^{-1/2}BA^{-1/2})^{q}  A^{1/2}. \nonumber
$$
\noindent The path $f_q(A,B)$ is the generalized geometric mean of $A$ and $B.$

{\bigskip}

In the present paper, we also study the representation of the generalized
 operator entropy which is defined for two positive operators A and B
on a Hilbert space and any real number $q \in ]0,1],$  by
\begin{equation}
S_{q}(A| B) = A^{1/2}
(A^{-1/2}BA^{-1/2})^{q} \ln (A^{-1/2}BA^{-1/2}) A^{1/2}. \nonumber
\end{equation}
\noindent The Bregman operator divergence introduced by Petz $[15]$ is defined by
\begin{equation}
D(A| B) = B-A-S(A| B).  \nonumber
\end{equation}

\noindent In $[6],$ Isa et al. have generalized $D(A| B) $ as follows
\begin{equation}
D_{q}(A| B) = f_{-q} (A| B) - f_q(A| B)- S_{q}(A| B). \nonumber
\end{equation}

 \noindent At the end of our paper, we also express the continued fraction representation of the divergence operator $D_{q}(A| B).$
 \\
\noindent For simplicity and clearness, we restrict ourselves to positive definite matrices, but our results can be,
without special difficulties, projected to
the case of positive definite operators from an infinite dimensional Hilbert space into itself.\\
\\
The computation of $ S_{q}(A| B)$ and $D_{q}(A| B) $ from the
original definitions impose many difficulties by
virtue of the appearance of the rational exponents of the
matrices. One fundamental of this paper is to remove this
difficulty and reveal a practical method, by matrix
continued fraction.

{\bigskip}

\noindent\large{\textbf{2. Definitions and notations }}
\\
\\
\noindent Let $\mathcal{M}_{m}$ be the algebra of real square matrices, where $m$ is a positive integer $\geq 2.$

The functions of matrix arguments play a widespreased role in science and
engineering, with applications areas from nuclear magnetic resonance $[1].$
\\
\noindent Let $A \in \mathcal{M}_{m},$ for a function $f$ with a series representation on an open
interval containing the eigenvalues of $A$, we are able to define the matrix
function $f\left( A\right) $ via the Taylor series for $f,$ see $[6].$%
\\
\noindent For any $A,B\in \mathcal{M}_{m}$ with $B$ invertible, we write ${\dfrac{A}{B}%
}=B^{-1}A.$  \\
\noindent It is easy to verify that for all matrices C and D with C invertible we have
\begin{equation*}
\displaystyle{\frac{A}{B}}=\displaystyle{\frac{CA}{CB}} \neq
\displaystyle{\frac{AC}{BC}}.
\end{equation*}
\\
\noindent{\bf Definition 2.1. \ }
Let $\left( A_{n}\right) _{n\geq 0}$, $\left( B_{n}\right) _{n\geq 0}$ be
two nonzero sequences of $\mathcal{M}_{m}$. The continued fraction of $%
\left( A_{n}\right) $ and $\left( B_{n}\right) $ denoted by $K(B_{n}/A_{n})$
is the quantity\rule[-10pt]{0pt}{10pt}\newline
\begin{equation*}
A_{0}+\dfrac{B_{1}}{A_{1}+\dfrac{B_{2}}{A_{2}+\cdot \cdot \cdot }}=\left[
A_{0};\frac{B_{1}}{A_{1}},\frac{B_{2}}{A_{2}},\cdot \cdot \cdot \right] .
\end{equation*}

\noindent Sometimes, we use briefly the notation $\left[ A_{0};\dfrac{B_{n}}{A_{n}}%
\right] _{n=1}^{+\infty }$. \\
\noindent The fractions $\dfrac{B_{n}}{A_{n}}$ and $\dfrac{P_{n}}{Q_{n}}=\left[ A_{0};%
\dfrac{B_{k}}{A_{k}}\right] _{k=1}^{n}$ are called, respectively, the $%
n^{th} $ partial quotient and the $n^{th}$ convergent of the continued
fraction $K(B_{n}/A_{n}).$
\\
\noindent Let $I$ be the $m^{th}$ order identity matrix, when $B_{n}=I$ for all $n\geq 1$, then $K(I/A_{n})$ is called a simple
continued fraction.
\\
\\
\noindent The continued fraction $K(B_{n}/A_{n})$ converges in ${\mathcal{M}}_{m}$ if
the sequence \\
$(F_{n})=({\frac{P_{n}}{Q_{n}}})$ converges
in ${\mathcal{M}}_{m}$ in the sense that there exists a matrix $F\in {
\mathcal{M}}_{m}$ such that $\lim_{n\rightarrow +\infty }||F_{n}-F||=0.$ In this case, we note%
\begin{equation*}
F=\left[ A_{0};\dfrac{B_{n}}{A_{n}}\right] _{n=1}^{+\infty }.
\end{equation*}

\noindent The following proposition gives an adequate method to calculate $K(B_{n}/A_{n}).$
\\
\\
\noindent {\bf Proposition 2.1. \ }
The elements $\left( P_{n}\right) _{n\geq -1}$ and $\left( Q_{n}\right)
_{n\geq -1\text{ }}$of the $n^{th}$ convergent of $K(B_{n}/A_{n})$ are given
by the relationships%
\begin{equation*}
 \left\{
\begin{array}{c}
P_{-1}=I,\quad P_{0}=A_{0} \\
Q_{-1}=0,\quad Q_{0}=I%
\end{array}%
\right. \text{ \ and \ \ }\left\{
\begin{array}{c}
P_{n}=A_{n}P_{n-1}+B_{n}P_{n-2} \\
Q_{n}=A_{n}Q_{n-1}+B_{n}Q_{n-2},%
\end{array}%
\right.  \ \ \ \ n\geq 1.
\end{equation*}%
\noindent{ Proof. } We prove it by induction.
\\

\noindent
The proof of the next Proposition is elementary and we left it to the reader.%

{\bigskip}

\noindent {\bf Proposition 2.2. \ }
(i) For any two matrices $C$ and $D$ with $C$ invertible, we have
\begin{equation*}
C\left[ A_{0};{\frac{B_{k}}{A_{k}}}\right] _{k=1}^{n}D=\left[ CA_{0}D  \
;{\frac{B_{1}D}{A_{1}C^{-1}},\frac{B_{2}C^{-1}}{A_{2}}, \frac{B_{k}}{A_{k}}}
\right]_{k=3}^{n}.
\end{equation*}
\noindent (ii) If two matrices $A $ and $B$ are similar, with $A=XBX^{-1},$ then
$f(A)$ and $f(B)$ are similar and we have
$f(A)=Xf(B)X^{-1}.$
\\
\\
\noindent In order to simplify the statement on some partial quotients of continued
fractions with matrix arguments, we need the following proposition which
is an example of equivalent continued fractions.\rule[-10pt]{0pt}{10pt}%
\newline

\noindent {\bf Proposition 2.3. \ }  Let \ $\left[ A_{0};\dfrac{B_{k}}{A_{k}}\right] _{k=1}^{+\infty}$ be
a given continued fraction. Then we have%
\begin{equation*}
\frac{P_{n}}{Q_{n}}=\left[ A_{0};\frac{B_{k}}{A_{k}}\right] _{k=1}^{n}=%
\left[ A_{0};\frac{X_{k}B_{k}X_{k-2}^{-1}}{X_{k}A_{k}X_{k-1}^{-1}}\right]
_{k=1}^{n}
\end{equation*}

\noindent where $X_{-1}=X_{0}=I$ and $X_{1}, X_{2},...X_{n}$ are arbitrary invertible matrices.
\\
\noindent{\bf Proof. \ }
Let $\dfrac{P_{n}}{Q_{n}}$ and $\dfrac{\widetilde{P}_{n}}{\widetilde{Q}_{n}}$%
 be the $n^{th}$ convergent of the continued fractions $\left[ A_{0};%
\dfrac{B_{k}}{A_{k}}\right] _{k=1}^{+\infty }$ and $\left[ A_{0};\dfrac{%
X_{k}B_{k}X_{k-2}^{-1}}{X_{k}A_{k}X_{k-1}^{-1}}\right] _{k=1}^{+\infty }$
respectively. By proposition $2$, for all $n\geq 1,$ we can write%
\begin{equation*}
\widetilde{P}_{n}=X_{n}A_{n}X_{n-1}^{-1}\widetilde{P}%
_{n-1}+X_{n}B_{n}X_{n-2}^{-1}\widetilde{P}_{n-2},
\end{equation*}

\noindent which is equivalent to%
\begin{equation*}
X_{n}^{-1}\widetilde{P}_{n}=A_{n}\left( A_{n}X_{n-1}^{-1} \widetilde{P}%
_{n-1}\right) +B_{n}( X_{n-2}^{-1} \widetilde{P}_{n-2}) .
\end{equation*}

\noindent This last result joined to the initial conditions prove that for all $n\geq
1,$ $X_{n-1}^{-1}\widetilde{P}_{n}=P_{n}.$ A similar result can be obtained for $Q_{n}$. Consequently, both continued
fractions have the same convergent and Proposition $3$ follows.

{\bigskip}

\noindent {\bf Definition 2.3. \ }
 Let $A$ be a positive definite matrix in ${\mathcal{M}}_{m}, X \in {\mathcal{M}}_{m}$ and $\alpha $ a
real number such that $0<\alpha <1.$ We define the matrix $A^{\alpha }$ by the formulae
\begin{equation*}
A^{\alpha }=\exp ( \alpha \ln A),
\end{equation*}
\noindent  where $"\exp"$ is the matrix exponential given by the series
\begin{equation*}
\exp(Y) =\sum_{n=0}^{+\infty }\frac{Y^{n}}{n!}
\end{equation*}

\noindent and $\ln $ is the neperian logarithm defined by

\begin{equation*}
\ln (A)=-2\sum_{n=0}^{+\infty }\frac{1}{2n+1}( \frac{I-A}{I+A})
^{2n+1}.
\end{equation*}

{\bigskip}

\noindent {\large{\bf 3. Main results }}

{\bigskip}

\noindent {\large{\bf 3.1 Continued fraction representation of $A^q \ln(A)$ }}
\\
\\
\noindent  This section is devoted to give a continued fraction representation of $A^{q } \ln(A),$ where $A$ is a positive definite matrix.

{\bigskip}

\noindent {\bf Theorem 3.1 \ }
 {\it Let $A\in {\mathcal{M}}_{m}$ be a positive definite matrix and $q$
a positive real number such that $ 0< q <1$. If we put
$A^{q }= \left[I; \dfrac{I}{A_k} \right]_{k=1}^{+\infty},  $
$ \ln A = \left[I; \dfrac{I}{A_k} \right]_{k=1}^{+\infty}$
 and $\varphi (A)
=\dfrac{I-A}{I+A}$ then, the continued fraction
expansions of
$ A^{q }\ln (A)$ is given by
\begin{equation*}
A^{q } \ln(A)=  \left[ 0;\dfrac{A_1+I}{A_1\widetilde{A}_{1} }
,\frac{ (A_1 A_2+I)(A_1+I) +\widetilde{A}_{1} \widetilde{A}_{2}}{(A_1A_2+A_2+I)\widetilde{A}_{2}},
\\
\\
\frac{(A_1\widetilde{A}_{1})^{-1}E_2F_3}{E_3-F_3}, \dfrac{ E_{n-1}F_{n}}{E_{n}-F_{n}}\right]_{n=4}^{+\infty }
\end{equation*}

\noindent where
\begin{eqnarray*}
 \left\{
 \begin{array}{c}
A_{1}=\frac{-I-q \varphi( A) }{2q \varphi( A) }, A_2= \dfrac{-6 qI}{(q^2-1)\varphi(A)},  \ \ \ \ \  \\
\\
  A_{2k}=\frac{-2q (q^{2}-2^{2})... ( q^{2}-( 2k-2)^{2}) }{( q^{2}-1)...
   (q^{2}-( 2k-1)^{2}) }( 4k-1) I, \ \ k\geq 2, \ \ \ \ \  \\
 \\
 A_{2k+1}=\dfrac{-( q^{2}-1) ... ( q^{2}-(
 2k-1) ^{2}) }{2q \left( q^{2}-2^{2}\right)...(q^{2}-4k^{2}) }( 4k+1) \varphi ^{2}( A), k\geq 1.
 \end{array}
 \right.
\end{eqnarray*}

\noindent For all $n \geq 1,$ the expression of $\widetilde{A}_{n}$ are given by the next relationships.

\begin{eqnarray*}
\left\{
 \begin{array}{c}
\widetilde{A}_{1}=(-2\varphi(A))^{-1}, \widetilde{A}_{2} = \dfrac{6I}{\varphi(A)} \\
  \\
\widetilde{A}_{2k}
=  \dfrac{-2(2.4... (2k-2))^2 }{(3.5...(2k-1))^2} (4k-1)I, \  k \geq 2  \\
  \\
\widetilde{A}_{2k+1}= \dfrac{(1.3.5... (2k-1) )^2}{( 2.4...(2k))^2 \varphi(A)}(4k+1)I,  \  k\geq 1.
\end{array}
\right.
\end{eqnarray*}

\noindent We also define
\begin{equation*}
 \left\{
 \begin{array}{c}
  E_{n}=Q_{n} \widetilde{Q}_{n}( Q_{n-2} P_{n-1} + \widetilde{Q}_{n-1} \widetilde{P}_{n-2}),\\
  F_{n}=Q_{n-2} \widetilde{Q}_{n-2} (
Q_{n} P_{n-1} + \widetilde{Q}_{n-1} \widetilde{ P}_{n}).
\end{array}
\right.
\end{equation*}

\noindent The matrices $P_n$ and $Q_{n}$ (resp. $ \widetilde{P}_{n}$ and $ \widetilde{Q}_{n}) $ are
the numerator and denominator of the $n^{th}$ convergent of $A^{q}$
(resp. $ \ln (A)$) which are defined for all $n\geq 1$ by
\begin{equation*}
 \left\{
 \begin{array}{c}
 P_{n} = A_{n}P_{n-1} +P_{n-2},
\\
Q_{n} =A_{n} Q_{n} + Q_{n-2},
\end{array}
\right. \ \text{  and }
\left\{
\begin{array}{c}
 \widetilde{ P}_{n} = \widetilde{A}_{n} \widetilde{P}_{n-1}  + \widetilde{ P}_{n-2}, \\
 \\
\widetilde{ Q}_n = \widetilde{A}_{n} \widetilde{Q}_{n-1}  + \widetilde{Q}_{n-2}.
\end{array}
\right.
\end{equation*}
}
\noindent In order to prove Theorem $3.1$, we begin by studying the real case.

{\bigskip}

\noindent{\large{\bf  3.2 The real case}}

{\bigskip}

\noindent We begin by giving some lemmas concerning the real continued
fraction which are important in the sequel.
The following lemma characterizes equivalence of continued fractions.

{\bigskip}

\noindent
{\bf Lemma 3.1 \ }$[7]$
{\it Let $(a_{n})_{n \geq 0}$ and $(b_{n})_{n \geq 1}$ be two non-zero sequences of real numbers.
The continued fractions

\begin{equation*}
\left[a_{0};\frac{b_{1}}{a_{1}}, \frac{b_{2}}{a_{2}},...
\frac{b_{n}}{a_{n}},...\right]=
\left[a_{0} \frac{1}{a_{1}^*},
\frac{1}{a_{2}^*},...,
\frac{1}{a_{n}^*},...\right]
\end{equation*}

\noindent where
 \begin{equation*}
\left\{
\begin{array}{c}
a_{1}^*={\frac{a_{1}}{b_{1}}}, \; a_{2}^*={\frac{b_1}{b_2}}a_2, \ \ \  \\
\\
a_{2k}^*={\frac{b_{1}b_{3} ... b_{2k-1}}{
b_{2}b_{4} ...  b_{2k}}}a_{2k},  \  k\geq 2, \\
 \\
a_{2k+1}^*={\frac{b_{2}b_{4}... b_{2k}}{
b_{1}b_{3}... b_{2k+1}}}a_{2k+1} , \  k\geq 1.
\end{array}
\right.
\end{equation*}
}
\noindent We now give a lemma which expresses the $n^{th}$ convergent for
the product of two continued fractions.
\\
\\
\noindent {\bf Lemma 3.2 } $[10] \ $ {\it  Let $C$ and $D$ be two real continued fractions which are defined by%
\begin{equation*}
C = \left[ c_{0}; \frac{1 }{ c_{1}}, \frac{1 }{ c_{n} }\right]_{n=2}^{+\infty}, \
D = \left[ d_{0}; \frac{1 }{ d_{1} }, \frac{1 }{ d_{n} }\right]_{n=2}^{+\infty},
\end{equation*}

\noindent where $c_{k}$ and $d_{k}$ are non-zero real numbers. If we put
\begin{equation*}
C_{n} = \left[ c_{0}; \frac{1 }{ c_{1} }, \frac{1 }{ c_{2} },...,\frac{1 }{ c_{n} } \right]
=\frac{^{c}p_{n}}{^{c}q_{n}} \ \text{ and } \ D_{n} = \left[ d_{0}; \frac{1 }{ d_{1} }, \frac{1 }{ d_{2} },...,\frac{1 }{ d_{n} }\right]=
\frac{^{d}p_{n}}{^{d}q_{n}},
\end{equation*}

\noindent then, for all $n\geq 1,$ we have
\begin{equation*}
C_n D_n =\left[ c_{0}d_{0};\frac{c_0c_{1}+d_0d_{1}+1}{%
c_{1}d_{1}},\frac{c_{1}d_{1}f_{2}}{e_{2}-f_{2}},\frac{e_{2}f_{3}}{e_{3}-f_{3}%
},...,\frac{e_{n-1}f_{n}}{e_{n}-f_{n}}\right] ,
\end{equation*}%

\noindent where
\begin{equation*}
\left\{
\begin{array}{c}
e_{n}=^{c}q_{n}^{d}q_{n}( ^{c}q_{n-2} \
^{c}q_{n-1}+ ^{d}q_{n-2} ^{d}q_{n-1}), \\
\\
f_{n}=^{c}q_{n-2}^{d}q_{n-2}( ^{c}q_{n-1} \ %
^{c}q_{n}+^{d}q_{n-1}^{d}q_{n}).%
\end{array}
\right.
\end{equation*}
}
\noindent The following Lemma gives two equivalent continued fraction expansions of $\lambda^{q}$,
where $\lambda$ and $q $ are two strictly positive real numbers.

{\bigskip}

\noindent {\bf Lemma 3.3.  \ } (i)  {\it Let $\lambda$ \ and $q$ be two positive real numbers, \newline
$ \varphi(\lambda) =\dfrac{1-\lambda}{1+\lambda}.$ The continued fraction expansions of $
\lambda^{q }$ is

\begin{equation*}
\lambda^{q }=\left[ 1;\frac{2q\varphi ( \lambda) }{-1-q
\varphi( \lambda) },\frac{( q^{2}-( k-1)
^{2}) \varphi ^{2}( \lambda) }{-( 2k-1) }\right]
_{k=2}^{+\infty }.
\end{equation*}

\noindent  (ii) The simple continued fraction of $\lambda^{q}$
is given by

\begin{equation*}
\lambda^{q }=\left[ 1;\frac{1}{c_{1}^*},\frac{1}{c_{2}^*},
... ,\frac{1}{c_{n}^*},... \right]
\end{equation*}
\noindent where
\begin{equation*}
\left\{
\begin{array}{c}
c_{1}^*={ \dfrac{-1-q \varphi( \lambda) }{2q
\varphi ( \lambda) }}, \ c_{2}^*={\dfrac{-6q}{(q^2-1)\varphi( \lambda) }}  \ \  \ \ \  \\
\\
c_{2k}^*
= {\dfrac{-2q \left( q^{2}-2^{2}\right)... \left( q^{2}-\left( 2k-2\right)
^{2}\right) }{\left(q ^{2}-1\right) ... \left(q ^{2}-\left( 2k-1\right)
^{2}\right) }}\left( 4k-1\right),  \  k\geq 2  \ \ \ \  \ \  \\
 \\
c_{2k+1}^*={\dfrac{-\left(q ^{2}-1\right) ...t \left( q ^{2}-\left( 2k-1\right)
^{2}\right) }{2q \left(q ^{2}-2^{2}\right) ... \left( q ^{2}-(2k)^{2}\right) }}\left(
4k+1\right) \varphi ^{2}\left( \lambda \right),\  k\geq 1 \ \ \ \  \ \
\end{array}
\right.
\end{equation*}
}
\noindent{Proof. } \ (i) See $[11].$
\\
\noindent (ii) By appropriate iteration and by applying Lemma $3.1$  we prove it.
\\
\\
\noindent {\bf Lemma 3.4 \ }
(i) {\it Let $\lambda $ be a real number such that $ \lambda
>0,$  $ \lambda \neq 1$ and $ \varphi(\lambda) =\frac{1-\lambda}{1+\lambda}.$ A continued fraction expansion of $\ln (\lambda)$ is given by
\begin{equation*}
\ln (\lambda) = \left[ 0; {- \frac{2\varphi(\lambda)
}{1}},  {-\frac{\varphi(\lambda)^2}{3}},\displaystyle{\frac{-2^2\varphi(\lambda)^2}{5}},
\displaystyle{\frac{-n^2
\varphi(\lambda)^2}{2n+1}}\right]_{n=3}^{+\infty}.
\end{equation*}

\noindent (ii) The simple continued fraction of $\ln (\lambda)$
is
\begin{equation*}
\ln (\lambda)=\left[ 0 ;\frac{1}{d_{1}^*},\frac{1}{d_{2}^*},
... ,\frac{1}{d_{n}^*},... \right]
\end{equation*}

\noindent where
\begin{equation*}
\left\{
\begin{array}{c}
d_{1}^*={\dfrac{1 }{-2 \varphi(\lambda)}},  \ \  d_{2}^*={\dfrac{6}{\varphi( \lambda) }},  \ \ \ \ \  \ \ \ \ \\
\\
d_{2k}^*
= {\dfrac{-2 ( 2.4..t (2k-2))^2 }
{(1.3.5...(2k-1) )^2}} (4k-1),  \  k\geq 2,  \ \ \ \ \ \ \   \\
 \\
d_{2k+1}^*={\dfrac{(1.3.5... (2k-1) )^2}{( 2.4... (2k))^2 \varphi(\lambda) }}(4k+1), \ k\geq 1. \ \ \ \ \ \
\end{array}
\right.
\end{equation*}.
}

\noindent { Proof. } \ (i) See $[12].$
\\
\noindent (ii) We deduce it by applying Lemma $3.1.$

{\bigskip}

\noindent The next Theorem is a real version of the previous Theorem $3.1.$

 {\bigskip}

\noindent {\textbf {Theorem 3.2 \ } {\it
With the same notations as bellow, let $\lambda$ and $q$ be two strictly positive real numbers such that
 $0<q<1.$ A continued fraction
representation of the real $\lambda^{q }\ln (\lambda)$\ is given by:%
\begin{equation*}
\lambda ^{q } \ \ln (\lambda)=\left[ 0; {\frac{{c_{1}^{* }+1}}{c_{1}^{* } {d}_{1}^{* }}},
{\frac{(c_{1}^{* } c_{2}^{* }+1)(c_{1}^{* }+1)+{d}_{1}^{* }d_{2}^{* }}
{(c_{1}^{* } c_{2}^{* }+c_{2}^{* }+1){d}_{2}^{* }}},
{\frac{(c_{1}^{* }{d}_{1}^{* })^{-1}e_2f_3}{e_3-f_3}},...,
{\frac{e_{n-1} f_{n}}{e_{n}-f_{n}}}\right]_{n= 4}^{+\infty},
\end{equation*}

\noindent where
\begin{equation*}
\left\{
\begin{array}{c}
e_{n} =q_{n} \widetilde{q}_{n} \left(
 {q}_{n-2}p_{n-1}+\widetilde{q}_{n-1} \widetilde{p}_{n-2} \right)
\\
f_{n} =q_{n-2} \widetilde{q}_{n-2} \left(
q_{n} p_{n-1} +\widetilde{q}_{n-1} \widetilde{p}_{n} \right),%
\end{array}
\right.
\end{equation*}
\\
\\
\noindent $p_n $ and $q_{n}$  (resp. $\widetilde{p}_{n})$ and $\widetilde{q}_{n})$ are numerator and denominator of the $n^{th}$ convergent of $\lambda
^{q } $ (resp. $\ln (\lambda)). $ They are
defined by
\begin{equation*}
\left\{
\begin{array}{c}
 p_{n} =c_{n}^{\ast }p_{n-1}
+p_{n-2}    \ \ \ \ \ \ \ \ \  \ \ \ \            \\
q_{n}=c_{n}^{\ast }q_{n-1}
+q_{n-2}   \ \ \ \ \ \ \ \ \ \ \ \  \
\end{array}
, \ \ \ \left\{
\begin{array}{c}
\widetilde{p}_{n}={d}_{n}^{\ast }\widetilde{p}_{n-1} +\widetilde{p}_{n-2} \ \ \ \ \ \ \ \ \ \ \  \\
\widetilde{q}_{n}={d}_{n}^{\ast }\widetilde{q}_{n-1}+ \widetilde{q}_{n-2}. \ \ \ \ \ \ \  \  \ \
\end{array}%
\right.
 \right.
\end{equation*}
}
\noindent {\bf Proof. \ } We apply Lemmas $3.2,$ $3.3$ and $3.4$ by putting $C=\lambda^q$ and $D=\ln(\lambda).$

{\bigskip}

\noindent {\large{\bf 3.3 Proof of Theorem $3.1.$ \ }}

{\bigskip}

\noindent Let $A \in {\mathcal{M}}_{m}$ be a positive definite matrix. Then there exists an invertible matrix
 $X$ such that
$A=XDX^{-1}$ where $D=diag( \lambda _{1},\lambda _{2},.., \lambda _{m})$   and $\lambda _{i}>0,$ for
$1\leq i<m.$
\\
\noindent As the functions $f(z)= z^q$ and $g(z)=\ln (z)$ are continuous in the open interval
$ \mathbb{R}_+^*, $ then we get
\begin{equation*}
A^{q } \ln(A)=X D^{q  } \ln(D) X^{-1}.
\end{equation*}

\noindent Let us define the sequences $\left( P_{n}\right) $ and $\left( Q_{n}\right) $
the numerator and denominator of the $n^{th}$ convergent of $D^{q }\ln (D)$ by
\begin{equation*}
\left\{
\begin{array}{c}
P_{1}=D_{1} + I,  P_{2}=( E'_{2}-F'_{2}) (D_1+I), P_3=(E'_3-F'_3)P_2+(D_1\widetilde{D}_{1})^{-1}E'_2F'_3P_1,  \  \  \ \ \
\\
Q_{1}=D_{1} \widetilde{D}_{1},  Q_{2}=( E'_{2}-F'_{2}) Q_{1}+D_{1} \widetilde{D}_{1} F'_{2},  Q_3=(E'_3-F'_3)Q_2+(D_1\widetilde{D}_{1})^{-1}E'_2F'_3Q_1   \  \
\end{array}%
\right.
\end{equation*}

\noindent and for $ n \geq 4,$
\begin{equation*}
\left\{
\begin{array}{c}
 P_{n}=\left( E'_{n}-F'_{n}\right) P_{n-1}+E'_{n-1} F'_{n} P_{n-2}
 \\
Q_{n}=\left( E'_{n}-F'_{n}\right) Q_{n-1}+E'_{n-1} F'_{n} Q_{n-2}%
\end{array}%
\right.
\end{equation*}
\begin{equation*}
\left\{
\begin{array}{c}
E'_{n}=Q_{n} \widetilde{Q}_{n} (
Q_{n-2}P_{n-1} +\widetilde{Q}_{n-1}\widetilde{P}_{n-2})
\\
F'_{n}=Q_{n-2} \widetilde{Q}_{n-2}\left(
Q_{n} P_{n-1} +\widetilde{Q}_{n-1} \widetilde{P}_{n} \right).%
\end{array}%
\right.
\end{equation*}
\\
\noindent The matrices $P_{n}, Q_{n}$ (resp. $\widetilde{P}_{n}, \widetilde{Q}_{n} )$ are the numerator and denominator of the $n^{th}$
convergent of $D^q$ (resp. $\ln(D) $) which are defined by

\begin{equation*}
\left\{
\begin{array}{c}
 P_{n} =D_{n}P_{n-1}
+P_{n-2}    \ \ \ \ \ \        \\
Q_{n}=D_{n}Q_{n-1}
+Q_{n-2}   \ \ \ \ \ \  
\end{array}%
\text{ \  and  \ }\left\{
\begin{array}{c}
\widetilde{P}_{n}=\widetilde{D}_{n}\widetilde{P}_{n-1} +\widetilde{P}_{n-2}
 \ \ \ \ \ \ \ \ \ \ \ \ \ \  \ \ \  \   \\
\widetilde{Q}_{n}=\widetilde{d}_{n}\widetilde{Q}_{n-1}+ \widetilde{Q}_{n-2}.
 \ \ \ \ \ \ \ \ \ \ \ \ \ \  \ \ \  \
\end{array}%
\right. \right.
\end{equation*}

\noindent We recall that

\begin{equation*}
\left\{
\begin{array}{c}
D_{1}=\dfrac{-I-q \varphi ( D) }{2q \varphi (
D) }, D_{2}=\dfrac{6qI}{(q^2-1)\varphi(D)}   \ \  \ \ \ \ \  \\
D_{2k}=\dfrac{-2q ( q^{2}-2^{2}) ... (q^{2}-( 2k-2)
^{2}) }{( q^{2}-1)... ( q^{2}-( 2k-1)
^{2}) }( 4k-1) I, \ k\geq 1 \ \ \ \ \ \ \ \ \
\\
D_{2k+1}=\dfrac{-( q^{2}-1) ... ( q^{2}-( 2k-1)
^{2}) }{2q( q^{2}-2^{2}) ... ( q^{2}-4k^{2}) }(
4k+1) \varphi ^{2}( D), \ k\geq 1. \ \
\end{array}%
\right.
\end{equation*}

\noindent and
\begin{eqnarray*}
\left\{
\begin{array}{c}
\widetilde{D}_{1}=(-2 \varphi(D))^{-1},  \widetilde{D}_{2} ={\frac{6I}{\varphi( D) }},  \ \ \  \ \ \ \ \  \ \ \\
\\
\widetilde{D}_{2k}
= {\dfrac{-2 ( 2.4...t (2k-2))^2 }
{(1.3.5...(2k-1) )^2}} (4k-1)I, \  k\geq 2  \ \ \ \ \ \ \   \\
 \\
\widetilde{D}_{2k+1}={\dfrac{(1.3.5... (2k-1))^2}{( 2.4...(2k))^2 \varphi(D) }}(4k+1)I,   \ k\geq 1. \ \ \ \ \ \ \
\end{array}
\right.
\end{eqnarray*}
\noindent We see that $P_{n},$  $Q_{n},$  $E'_n$ and $F'_n$ are diagonal matrices, so
we put%
\begin{equation*}
\left\{
\begin{array}{c}
P_{n} = diag( p_{n}^{1}, \ p_{n}^{2},..., p_{n}^{m}) \\
Q_{n} =diag( q_{n}^{1}, \  q_{n}^{2},
..., q_{n}^{m})
\end{array}
, \  \ \left\{
\begin{array}{c}
E'_{n} = diag( e_{n}^{1}, \ e_{n}^{2},
..., e_{n}^{m}),
 \\
F'_{n} = diag( f_{n}^{1}, \ f_{n}^{2},
..., f_{n}^{m}).
\end{array}
\right. \right.
\end{equation*}

\noindent We obtain for each $1\leq i\leq m,$
\begin{equation*}
\left\{
\begin{array}{c}
p_{1}^{i}=c_{i1}^{\ast } + 1, \text{ \ } p_{2}^{i}= (e_{2}^i-f_{2}^i) (c_{i1}^{\ast }+1),
p_3^i=(e_3^i-f_3^i)p_2+(c_1^{\ast}{d}_{1}^{\ast})^{-1}p_1^i \text{ \  \  } \\
\\
q_{1}^{i}= c_{i1}^{\ast }
{d}_{i1}^{\ast }, \  \  q_{2}^{i}= (e_{2}^i-f_{2}^i)q_1^i+c_{i1}^{\ast }d_{i1}^{\ast} f_2^i, q_3^i=(e_3^i-f_3^i)q_2+(c_1^{\ast}{d}_{1}^{\ast})^{-1}q_1^i %
\end{array}%
\right.
\end{equation*}

\noindent and for $n \geq 4,$ we have
\rule[-10pt]{0pt}{10pt}
\begin{equation*}
\left\{
\begin{array}{c}
p_{n}^{i}=( e_{n}^{i }-f_{n}^{i })
p_{n-1}^{i}+e_{n-1}^{i } f_{n}^{i } p_{n-2}^{i}, \\
\\
q_{n}^{i}=( e_{n}^{i }-f_{n}^{i })
q_{n-1}^{i}+e_{n-1}^{i} f_{n}^{i } q_{n-2}^{i},%
\end{array}%
\right.
\end{equation*}%
\\
\begin{equation*}
\left\{
\begin{array}{c}
e_{n}^{i }=q_{n}^{i} \widetilde{q}_{n}^{i}( q_{n-2}^{i} p_{n-1}^{i}+ \widetilde{q}_{n-1}^{i} \widetilde{p}_{n-2}^{i}) \\
{\medskip } \\
f_{n}^{i }=q_{n-2}^{i} \widetilde{q}_{n-2}^{i}( q_{n}^{i} p_{n-1}^{i} + \widetilde{q}_{n-1}^{i} \widetilde{p}_{n}^{i} ). %
\end{array}%
\right.
\end{equation*}
\\
\noindent By Lemma $3.2,$ we deduce that $\dfrac{p_{n}^{i}}{q_{n}^{i}}$ converges to $ \lambda _{i}^q \ln(\lambda_i)$ for $1\leq i\leq m.$
It follows that the matrix $\dfrac{P_{n}}{Q_{n}}$ converges to $D^{q} \ln (D).$ So, we get

\begin{equation*}
D^{q} \ln(D)=\left[ 0;\frac{D_{1}+I}{D_{1}\widetilde{D}_{1}},
\dfrac{(D_{1}\widetilde{D}_{1})^{-1}E'_{2}F'_3}{E'_{3}-F'_{3}},
\frac{%
E'_{n-1}F'_{n}}{E'_{n}-F'_{n}}\right] _{n= 4}^{+\infty}.
\end{equation*}
\\
\noindent By Proposition $3.2,$ we have%
\begin{eqnarray*}
A^{q} \ln(A) &=&X\left( D^{q } \ln(D) \right) X^{-1}
\\
&=&\left[0 ;\frac{(D_{1}+I)X^{-1}}{D_{1}\widetilde{D}_1X^{-1}},\frac{(D_{1}\widetilde{D}_{1})F'_{2}X^{-1}}{E'_{2}-F'_{2}},
\dfrac{(D_{1}\widetilde{D}_{1})^{-1}E'_{2}F'_3}{E'_{3}-F'_{3}},
\frac{ E'_{n-1}F'_{n}}{E'_{n}-F'_{n}}\right] _{n=4}^{+\infty}
\end{eqnarray*}

\noindent Let us define the sequence $(X_n)_{n\geq -1}$ by
$X_{-1}=X_{0}=I,$ and for all $n\geq 1,$ $X_{n}=X.$ Then we have
\rule[-10pt]{0pt}{10pt}\newline
\begin{equation*}
\left\{
\begin{array}{c}
\dfrac{X_{1}(D_1+ I)X^{-1}X_{-1}^{-1}}{X_{1}D_1\widetilde{D}_1 X^{-1}X_{0}^{-1}}=\dfrac{A_1
+I}{A_1\widetilde{A}_1}, \\
\\
\dfrac{X_{2}(D_1\widetilde{D}_1F'_2X^{-1})X_{0}^{-1}}{X_{2}(E'_2-F'_2)X_{1}^{-1}}=
\dfrac{A_1\widetilde{A}_1F_2}{E_2-F_2}, \\
\\
\dfrac{X_{n}(E'_{n-1}F'_n)X_{n-2}^{-1}}{X_{n}(E'_n-F'_n)X_{n-1}^{-1}}=
  \dfrac{E_{n-1}F_n}{E_n-F_n}
\end{array},
\right.
\end{equation*}

\noindent with $XE'_{n}X^{-1}=E_n$ and $XF'_{n}X^{-1}=F_n$ for all $ n \geq 2.$
\\
\noindent By applying the result of Proposition $2.3$ to the sequence $
(X_{n})_{n\geq -1},$ we finish the proof of Theorem $3.1.$

{\bigskip}

\noindent {\large{\bf{3.4 Representation of the generalized operator entropy}}}

{\bigskip}

\noindent {\bf Theorem 3.3 } {\it  Let A and B be two positive definite matrices in ${\mathcal M}_m, $ $q$
a positive real number such that $0< q <1.$    A
continued fraction representation of the generalized operator
 entropy $S_q(A|B)$ is given by

 $$
  S_q(A|B)=
  $$
\begin{equation*}
\left[ 0; \displaystyle{
 \frac{(A'_1+I)A^{1/2}}{A'_1\widetilde{A'}_{1}A^{-1/2} }
,\frac{ ((A'_1 A_2+I)(A'_1+I) +\widetilde{A'}_{1}\widetilde{A'}_{2})A^{-1/2}}{(A'_1A'_2+A'_2+I)\widetilde{A'}_{2}},
\frac{(A'_1\widetilde{A'}_{1})^{-1}E_2F_3}{E_3-F_3}, \frac{ E_{n-1}F_{n}}{E_{n}-F_{n}}} \right]_{n=4}^{+\infty }
 \end{equation*}

\noindent where $ A'_k=A_k $ and $ \widetilde{A'}_{k}=\widetilde{A}_{k} $  which are defined in Theorem $3.1$ by
the equalities $(0.1)$ and $(0.2)$ by replacing
$\varphi(A) $ by
 $$ \varphi(A,B)=\varphi(A^{-1/2}BA^{-1/2})=A^{1/2}(\frac{A-B}{A+B})A^{-1/2}.
 $$
}
\noindent {\bf Proof. \ } We have
 $$ (A^{-1/2}BA^{-1/2})^q \ln (A^{-1/2}BA^{-1/2})=A^{-1/2} S_q(A|B)A^{-1/2}.$$

\noindent By applying Theorem $3.1,$ the continued fraction representation of $S_q(A|B)$ is

$$
S_q(A|B)=
$$
$$
A^{1/2} \left[ 0; \frac{A'_1+I}{A'_1\widetilde{A'}_{1} }
,\frac{ (A'_1 A_2+I)(A'_1+I) +\widetilde{A'}_{1}\widetilde{A'}_{2}}{(A'_1A'_2+A'_2+I)\widetilde{A'}_{2}},
\frac{(A'_1\widetilde{A'}_{1})^{-1}E_2F_3}{E_3-F_3}, \frac{ E_{n-1}F_{n}}{E_{n}-F_{n}}\right]_{n=4}^{+\infty }A^{1/2}.
$$
\noindent That is
$$
S_q(A|B)=
$$
$$
 \left[ 0; \frac{(A'_1+I)A^{1/2}}{A'_1\widetilde{A'}_{1}A^{-1/2} }
,\frac{( (A'_1 A_2+I)(A'_1+I) +\widetilde{A'}_{1}\widetilde{A'}_{2})A^{-1/2}}{(A'_1A'_2+A'_2+I)\widetilde{A'}_{2}},
\frac{(A'_1\widetilde{A'}_{1})^{-1}E_2F_3}{E_3-F_3}, \frac{ E_{n-1}F_{n}}{E_{n}-F_{n}}\right]_{n=4}^{+\infty }
$$

\noindent which completes the proof of Theorem $3.3.$

{\bigskip}

\noindent Let $q \in \mathbb{R}, $ such that $0<q<1$ and $n $ be an integer, we gave some properties of $S_q(A|B),$ which are shown in $[6].$
\\
\\
\noindent {\bf Lemma 3.5. \ }
 {\it Let $A $ and $B$ be two positive definite matrices in ${\mathcal M}_m, $ $ n \in \mathbb{N}^*,$ then we have

 $$ (i) S_n(A|B)=(BA^{-1})^n S(A|B)= S(A|B) (A^{-1}B)^{n}.$$
$$ (ii)
 S_{2n}(A|B)=(BA^{-1})^n S(A|B)(A^{-1}B)^{n}.$$
$$  (iii) S_{2n+1}(A|B)=(BA^{-1})^n S_1(A|B)(A^{-1}B)^{n}.$$
}
\noindent  Now we have the next result  which gives a continued fraction expansions of the relative operator
 entropy $S_n(A|B).$ \\
\\
\\
\noindent {\bf Corollary 3.4 } {\it Let A and B be two positive definite matrices in ${\mathcal M}_m, $ $n \in \mathbb{N}.$  A
continued fraction representation of the generalized operator
 entropy $S_n(A|B)$ is given by
\begin{equation*}
 S_n(A|B)=\left[ 0;
 \displaystyle{ \frac{2(BA^{-1})^nA \left( \frac{A-B}{A+B}\right)}{I}},
 \displaystyle{ \frac{-k^2A
\left(\frac{A-B}{A+B}\right)^2A^{-1}}{(2k+1)I}}
\right]_{k=1}^{+\infty}.
\end{equation*}
}
\noindent {\bf Proof. \ } In order to prove this Theorem, we recall a continued fraction expansion of
the relative operator entropy $S(A|B)$ (see $[14]$).
\begin{equation*}
 S(A|B)=\left[ 0;
 \displaystyle{ \frac{2A \left( \frac{A-B}{A+B}\right)}{I}},
\displaystyle{ \frac{-k^2A
\left(\frac{A-B}{A+B}\right)^2A^{-1}}{(2k+1)I}}
\right]_{k=1}^{+\infty}.
\end{equation*}

\noindent By applying Lemma $3.5,$ we complete the proof of Corollary $3.4.$

{\bigskip}

\noindent {\large{\textbf{4. The operator divergence }}}

{\bigskip}

\noindent The Bregman operator divergence is $D(A|B)=B-A-S(A|B).$  \\
We see that its positivity is
assured by
$$
S(A|B)=A^{1/2}\ln (A^{-1/2}BA^{-1/2}) A^{1/2} \leq A^{1/2}(A^{-1/2}BA^{-1/2}-I) A^{1/2}=B-A.
$$

\noindent In Isa et al. $[7]$ have generalized $D(A| B) $ as follows
\begin{equation}
D_{q}(A| B) = f_{-q} (A| B) - f_q(A| B)- S_{q}(A| B). \nonumber
\end{equation}

\noindent where
$$
 f_{q} (A| B)=  A^{1/2}( A^{-1/2}BA^{-1/2})^q A^{1/2}.
 $$
 \\
\newpage
\noindent {\it {\bf Theorem 4.1 } Let A and B be two positive definite matrices in ${\mathcal M}_m, $ $q$
a positive real number such that $0< q <1.$    A
continued fraction representation of the operator
 divergence $D_q(A|B)$ is given by
 \begin{eqnarray*}
 D_q(A|B)&=& \left[ A^{-1/2}BA^{-1/2}-I;
  \dfrac{(A'_1+I+\widetilde{A'}_0 \widetilde{A'}_1)A^{1/2}}{A'_1\widetilde{A'}_{1}A^{-1/2} },
  \right. \\
&& \left.
\frac{ ((A'_1 A'_2+I)(A'_1+I) +\widetilde{A'}_{1}\widetilde{A'}_{2})A^{-1/2}}{(A'_1A'_2+A'_2+I)\widetilde{A'}_{2}},%
\frac{(A'_1\widetilde{A'}_{1})^{-1}E_2F_3}{E_3-F_3}, \frac{ E_{n-1}F_{n}}{E_{n}%
-F_{n}}\right]_{n=4}^{+\infty}
 \end{eqnarray*}
\noindent  For $k \geq 2,$ $\widetilde{A'}_k$ are the same as in Theorem $3.3.$ The matrices
$A'_k, E_k$ and $F_k$ for all $k \geq 1$ are also the same as in Theorem $3.3.$
}
{\bigskip}

\noindent{\bf Proof. } Since $f_{-q}(A|B)=f_{1+q}(A|B),$ then we have

\begin{eqnarray*}
 D_q(A|B)&=& A^{1/2}((A^{-1/2}BA^{-1/2})^{q+1}- (A^{-1/2}BA^{-1/2})^q - \\
          & & (A^{-1/2}BA^{-1/2})^q \ln (A^{-1/2}BA^{-1/2}))  A^{1/2} \\
      &=& A^{1/2}(A^{-1/2}BA^{-1/2})^{q}(- \ln (A^{-1/2}BA^{-1/2})+ A^{-1/2}BA^{-1/2}-I) A^{1/2}
\end{eqnarray*}

\noindent So, the proof of Theorem $4.1$ is shown by the similar way of Theorem $3.3$  by replacing
$\ln (A^{-1/2}BA^{-1/2})$ by
$(-\ln (A^{-1/2}BA^{-1/2})+A^{-1/2}BA^{-1/2}-I).$

{\bigskip}

\noindent{\bf 5. Numerical applications }

{\bigskip}

This paragraph will provide some numerical data to illustrate the preceding
results. The focus will be on the results obtained for the generalized operator entropy and the operator divergence.

 {\bigskip}

\noindent {\bf 5.1. Numerical example of the $A^{q } \ln(A),$} \\
\\
 We start this section by giving an example to illustrate the theoretical results obtained in theorem $3.2.$

{\bigskip}

 \noindent{ \bf Exemple 1.}

 {\bigskip}

\begin{tabular}{|p{1pc}|p{1pc}|p{4pc}|p{4pc}|p{4pc}|p{4pc}|p{4pc}|}
\hline $\lambda$ & q   &$\lambda^{q}\ln\lambda-F_{1}$ &$\lambda^{q}\ln\lambda-F_{2}$ & $\lambda^{q}\ln\lambda-F_{3}$ & $\lambda^{q}\ln\lambda-F_{4}$ & $\lambda^{q}\ln\lambda-F_{5}$\\
\hline $2$ & $\frac{1}{2}$    &  $0.04692481$  &$0.1478 \; 10^{-2}$ &  $0.4462\; 10^{-4}$&$0.1331\; 10^{-5}$ & $0.394\; 10^{-7}$ \\
\hline $2$ & $\frac{1}{3}$    &  $0.03997739$  &$0.1269\; 10^{-2}$ &  $0.3845\; 10^{-4}$&$0.115\; 10^{-5}$ & $0.343\; 10^{-7}$ \\

\hline $3$ & $\frac{1}{2}$ &    $0.2361856$  &$0.01855478$ &  $0.1369\; 10^{-2}$&$0.9973\; 10^{-4}$ & $0.7224\; 10^{-5}$ \\
\hline
\end{tabular} \\

Now, we pass to a more general case than the previous one, matrix case. \\

{\bigskip}
\newpage
\noindent{ \bf Exemple 2.}

{\bigskip}

Let $A\in \mathcal{M}_{m}$ be a positive definite matrix, such that
$$A=\begin{pmatrix}
2 & 1 \\
1 & 3
\end{pmatrix}$$

  We calculate the difference between the exact value of $A^{q}lnA$ and its first convergent. We will take for example $q=\frac{1}{3},$  we get

  $$A^{q}lnA- F_{1}=\begin{pmatrix}
0.08867877674 & 0.1379741093 \\
0.1379741091 & 0.226652885
\end{pmatrix}$$

$$A^{q}lnA- F_{2}=\begin{pmatrix}
0.00934803094 & 0.0150871241 \\
0.0150871239 & 0.024435154
\end{pmatrix}$$

$$A^{q}lnA- F_{3}=\begin{pmatrix}
0.93618974\,10^{-3} & 0.15145295\,10^{-2} \\
0.15145293\,10^{-2} & 0.2450718\,10^{-2}
\end{pmatrix}$$

$$A^{q}lnA- F_{4}=\begin{pmatrix}
0.9212614\,10^{-4} & 0.1490607\,10^{-3} \\
0.1490605\,10^{-3} & 0.241186\,10^{-3}
\end{pmatrix}$$

$$A^{q}lnA- F_{5}=\begin{pmatrix}
0.900254\,10^{-5} & 0.145655\,10^{-4} \\
0.145653\,10^{-4} & 0.23567\,10^{-4}
\end{pmatrix}$$

{\bigskip}

\noindent {\bf 5.2. Numerical example of the generalized operator entropy }
\\
\\
Now, we turn to the first main objective of this work; approximation of the generalized operator entropy by using continued fractions. \\

{\bigskip}

\noindent{ \bf Exemple 3.}

{\bigskip}

Let $A$ and $B$ be two positive definite matrices in $\mathcal{M}_{m}$  such that \\

$$A=\begin{pmatrix}
5 & 0 & 0\\
0 & 2 & 1 \\
0 & 1 & 5
\end{pmatrix}  \text{ and }  B=\begin{pmatrix}
3 & 1&0 \\
1 & 4 & 1 \\
0 &1 & 3
\end{pmatrix}.$$

We calculate the difference between $S_{q}(A\mid B)$ and its first five convergents, we will take for example $q=\frac{1}{3},$ we got the following results \\

$$F_{1}-S_{q}(A\mid B)=\begin{pmatrix}
0.054190989 & -0.04186083842 & 0.0061362095 \\
-0.0418608386 & -0.09860205655 & -0.0089742420 \\
0.00613620937 & -0.0089742424 & 0.041304948
\end{pmatrix}$$

$$F_{2}-S_{q}(A\mid B)=\begin{pmatrix}
0.001377942 & -0.00151980232 & 0.0002817516 \\
-0.0015198025 & -0.00384514555 & -0.0004732079 \\
0.00028175148 & -0.0004732083 & 0.000786263
\end{pmatrix}$$

$$F_{3}-S_{q}(A\mid B)=\begin{pmatrix}
0.000032775 & -0.00005387752 & 0.86254\,10^{-5} \\
-0.0000538777 & -0.00014206055 & -0.0000193084 \\
0.862533\,10^{-5} & -0.0000193088 & 0.000014661
\end{pmatrix}$$

$$F_{4}-S_{q}(A\mid B)=\begin{pmatrix}
0.718\,10^{-6} & -0.190962\,10^{-5} & 0.2162\,10^{-6} \\
-0.19098\,10^{-5} & -0.516755\,10^{-5} & -0.7279\,10^{-6} \\
0.21615\,10^{-6} & -0.7283\,10^{-6} & 0.263\,10^{-6}
\end{pmatrix}$$

$$F_{5}-S_{q}(A\mid B)=\begin{pmatrix}
0.13\,10^{-7} & -0.6722\,10^{-7} & 0.65\,10^{-8} \\
-0.674\,10^{-7} & -0.18655\,10^{-6} & -0.272\,10^{-7} \\
0.643\,10^{-8} & -0.276\,10^{-7} & -0.2\,10^{-8}
\end{pmatrix}$$

{\bigskip}

\noindent{\bf 5.3 Numerical example of the operator divergence } \\
\\
We end this paragraph by illustrating the theoretical results concerning the divergence operator. \\

We keep the same data as we have in the previous example. We calculate the difference between $D_{q}(A\mid B)$ and its first five convergents, we got the following results

{\bigskip}

\noindent{\bf Exemple 4.}
\\
$$F_{1}-D_{q}(A\mid B)=\begin{pmatrix}
-0.069264959 & 0.03318154762 & -0.0099888506 \\
0.0331815487 & 0.07019678255 & 0.0044140547 \\
-0.00998885042 & 0.0044140552 & -0.048288369
\end{pmatrix}$$

$$F_{2}-D_{q}(A\mid B)=\begin{pmatrix}
-0.001821578 & 0.00117926062 & -0.0004227273 \\
0.0011792617 & 0.00278467055 & 0.0003187543 \\
-0.00042272712 & 0.0003187548 & -0.000933848
\end{pmatrix}$$

$$F_{3}-D_{q}(A\mid B)=\begin{pmatrix}
-0.000045885 & 0.00004116862 & -0.0000132690 \\
0.0000411697 & 0.00010356255 & 0.0000139272 \\
-0.00001326882 & 0.0000139277 & -0.000018017
\end{pmatrix}$$

$$F_{4}-D_{q}(A\mid B)=\begin{pmatrix}
-0.1115\,10^{-5} & 0.144362\,10^{-5} & -0.3629\,10^{-6} \\
0.14447\,10^{-5} & 0.378555\,10^{-5} & 0.5380\,10^{-6} \\
-0.36272\,10^{-6} & 0.5385\,10^{-6} & -0.350\,10^{-6}
\end{pmatrix}$$

$$F_{5}-D_{q}(A\mid B)=\begin{pmatrix}
-0.27\,10^{-7} & 0.5062\,10^{-7} & -0.100\,10^{-7} \\
0.507\,10^{-7} & 0.14255\,10^{-6} & 0.217\,10^{-7} \\
-0.982\,10^{-8} & 0.222\,10^{-7} & -0.3\,10^{-8}
\end{pmatrix}.
$$

{\bigskip}

\end{document}